\documentclass[12pt]{amsart}
\usepackage{amsfonts}
\usepackage{amsmath}

\newtheorem{theorem}{Theorem}
\theoremstyle{plain}
\newtheorem{definition}{Definition}
\newtheorem{lemma}{Lemma}

\newtheorem{remark}{Remark}

\numberwithin{equation}{section}

\begin{document}
\title{A pseudodifferential H\"{o}rmander's inequality}
\author{Chikh BOUZAR}
\address{Department of Mathematics, Oran-Essenia University, Algeria}
\email{bouzar@wissal.dz ; bouzar@yahoo.com}

\begin{abstract}
The classical H\"{o}rmander's inequality for linear partial differential
operators with constant coefficients is extended to pseudodifferential
operators.
\end{abstract}

\maketitle

\section{Introduction}

The classical H\"{o}rmander's inequality proved in \cite{H1} for linear
partial differential operators with constant complex coefficients $P(D)$\
and a bounded domain $\Omega $ of $\mathbb{R}^{n}\ $claims that
\begin{equation*}
\left\| \left( \partial ^{\alpha }P\right) \left( D\right) u\right\|
_{0}\leq c\left\| P\left( D\right) u\right\| _{0}\text{ , }\forall u\in
C_{0}^{\infty }\left( \Omega \right)
\end{equation*}
where $\alpha \in \mathbb{N}^{n}$ and\ $c>0$ is independent of $u.$

In the proof,\ we obviously see that in fact the following result hold
\begin{equation*}
\forall \delta >0,\exists \rho >0,\forall u\in C_{0}^{\infty }\left( \Omega
\right) ,diam(\Omega )<\rho ,
\end{equation*}
\begin{equation*}
\left\| \left( \partial ^{\alpha }P\right) \left( D\right) u\right\|
_{0}\leq \delta \left\| P\left( D\right) u\right\| _{0}
\end{equation*}

The aim of this paper is to give an extension of this H\"{o}rmander's
inequality to classical pseudodifferential operators with constant
coefficients. The proposed pseudodifferential H\"{o}rmander's inequality
includes the given cases in \cite{P} and \cite{B-K}. The formulation of our
result is the following theorem.

\begin{theorem}
Let $P\left( D\right) $ be a pseudodifferential operator of the class $S^{m},
$ and let $s\in \mathbb{R},\theta \geq 1$ and$\ \alpha \in \mathbb{N}^{n},$\
then
\begin{equation*}
\forall \delta >0,\exists \rho >0,\exists c>0,\forall u\in C_{0}^{\infty
}\left( \Omega \right) ,diam(\Omega )<\rho ,
\end{equation*}
\begin{equation*}
\left\| \left( \partial ^{\alpha }P\right) \left( D\right) u\right\|
_{s}\leq \delta \left\| P\left( D\right) u\right\| _{s}+c\left\| u\right\|
_{s+m-\theta }
\end{equation*}
$\ $
\end{theorem}

\section{Preliminaries}

The set $\Omega $\ denotes an open domain of $\mathbb{R}^{n}$\ and $\mathbb{N%
}$ the set of natural numbers $\left\{ 1,2,...\right\} .$\ The notations and
classical definitions from the theory of distributions and
pseudodifferential operators are given in \cite{H2} and \cite{T}; in
particular, $C_{0}^{\infty }(\Omega )$ is the space of infinitely
differentiable functions with compact support and $H^{s},s\in \mathbb{R},$\
is the Sobolev space on $\mathbb{R}^{n}$\ with scalar product and norm
denoted, respectively, $\left( \cdot ,\cdot \right) _{s}$\ and $\left\|
\cdot \right\| _{s}$ .

We will use the following classical inequalities, $s,t\in \mathbb{R},$%
\begin{equation*}
\left| \left( u,v\right) _{s}\right| \leq \left\| u\right\| _{s-t}\left\|
v\right\| _{s+t}\text{ , }\forall u\in H^{s-t},\forall v\in H^{s+t}\text{ \ ,%
}
\end{equation*}
and

\begin{lemma}
Let $\varphi \in C_{0}^{\infty }(\Omega )$ and $s\in \mathbb{R}$\ , then
there exists a constant $c>0$\ such that
\begin{equation*}
\left\| \varphi u\right\| _{s}\leq \max_{x}\left| \varphi (x)\right| \left\|
u\right\| _{s}+c\left\| u\right\| _{\sigma }\text{ , }\forall u\in H^{s}%
\text{\ ,}
\end{equation*}
where $\sigma <s-1.$
\end{lemma}

\begin{proof}
See lemma 2 of \cite{Pe}
\end{proof}

The class of symbols of pseudodifferential operators with constant
coefficients is defined as follows.

\begin{definition}
The class $S^{m}$ is the space of infinitely differentiable functions $P(\xi
)$ defined on $\mathbb{R}^{n}$\ and satisfying $\forall \alpha \in \mathbb{Z}%
_{+}^{n}$ , there exists $c>0$\ such that
\begin{equation*}
\left| \left( \partial ^{\alpha }P\right) \left( \xi \right) \right| \leq
c\left( 1+\left| \xi \right| \right) ^{m-\left| \alpha \right| }\text{ , }%
\forall \xi \in \mathbb{R}^{n}
\end{equation*}
\end{definition}

A pseudodifferential operator $P\left( D\right) $ of order $m$\ with
constant coefficients is an operator acting on functions $u\in H^{s}$\ by
the formula
\begin{equation*}
P\left( D\right) u=\int_{\mathbb{R}^{n}}e^{ix\cdot \xi }P\left( \xi \right)
\widehat{u}(\xi )d\xi
\end{equation*}
where $P(\xi)\in S^{m}$ is called the symbol of $P\left( D\right) ,$\ and $%
\widehat{u}$\ denotes the Fourier transform of $u$\ .\

The pseudodifferential operators $\left( \partial ^{\alpha }P\right) \left(
D\right) ,$ $\alpha \in \mathbb{N}^{n},$ and $\overline{P}\left( D\right) $
are the pseudodifferential operators with respective symbols $\left(
\partial ^{\alpha }P\right) (\xi )$\ and $\overline{P\left( \xi \right) }$.
It is easy to see that the $H^{s}$-adjoint operator of $P\left( D\right) $
is the operator $\overline{P}\left( D\right) ,$ and we have\
\begin{equation*}
\left\| \overline{P}u\right\| _{s}=\left\| Pu\right\| _{s}\text{ , }\forall
u\in H^{s}
\end{equation*}

\begin{remark}
Our a priori estimates are local and in view of a classical result of the
theory of pseudodifferential operators, all the operators considered in this
work are properly supported.
\end{remark}

By $\Omega _{\varepsilon }$\ we denote the open ball of center the origin
and radius $\varepsilon >0.$\ Let $\varphi \left( x\right) \in C_{0}^{\infty
}\left( \mathbb{R}^{n}\right) $ such that $\varphi \left( x\right) =1$ for $%
\left| x\right| \leq 1,0\leq \varphi \left( x\right) \leq 1$ and $\varphi
\left( x\right) =0$\ for $\left| x\right| >2.$ Define $\varphi _{\varepsilon
}\left( x\right) =\varphi (\frac{x}{\varepsilon })$ and let $P\left(
D\right) $\ be a pseudodifferential operator of the class\ $S^{m}$ . The
operator $[P\left( D\right) ,\varphi _{\varepsilon }]$\ denotes the
commutator of the pseudodifferential operator $P\left( D\right) $\ and the
operator of multiplication by the function $\varphi _{\varepsilon }\left(
x\right) .$

\begin{lemma}
If $0<\rho <\varepsilon ,$ then the operator $[P\left( D\right) ,\varphi
_{\varepsilon }]$ is of infinite order in $C_{0}^{\infty }(\Omega _{\rho }),$%
\ i.e. for every reals $s$ and $s^{\prime },$ there exists $c>0$\ such that
\begin{equation*}
\left\| \lbrack P\left( D\right) ,\varphi _{\varepsilon }]u\right\| _{s}\leq
c\left\| u\right\| _{s^{\prime }}\text{ \ , }\forall u\in C_{0}^{\infty
}(\Omega _{\rho })
\end{equation*}
\end{lemma}

\begin{proof}
It is deduced from the fact that the symbol of the pseudodifferential
operator $[P\left( D\right) ,\varphi _{\varepsilon }]$\ is identically
equals zero on a neighbourhood\ of the open set $\Omega _{\rho }$\ \
\end{proof}

\begin{remark}
We will apply the following algebraic inequality,
\begin{equation*}
2ab\leq \varepsilon a^{2}+\frac{1}{\varepsilon }b^{2}\text{ , }\forall
\varepsilon >0,\forall a\geq 0,\forall b\geq 0
\end{equation*}
\end{remark}

\section{The inequality}

The principal result of the paper is the following theorem which is an
extension of H\"{o}rmander's inequality.

\begin{theorem}
Let $P\left( D\right) $ be a pseudodifferential operator of the class $%
S^{m}, $ and let $s\in \mathbb{R},\theta \geq 1$ and$\ \alpha \in \mathbb{N}%
^{n},$\ then
\begin{equation*}
\ \forall \delta >0,\exists \rho >0,\exists c>0,\forall u\in C_{0}^{\infty
}\left( \Omega \right) ,diam(\Omega )<\rho ,
\end{equation*}
\begin{equation}
\left\| \left( \partial ^{\alpha }P\right) \left( D\right) u\right\|
_{s}\leq \delta \left\| P\left( D\right) u\right\| _{s}+c\left\| u\right\|
_{s+m-\theta }  \label{1}
\end{equation}
\end{theorem}

\begin{proof}
Without lost of generality let $\Omega _{\varepsilon }$\ be the
open ball of center the origin and radius $\varepsilon >0.$\ Let
$\varphi \left( x\right) \in C_{0}^{\infty
}(\mathbb{R}^{n}),\varphi \left( x\right) =1$ for $\left| x\right|
\leq 1,0\leq \varphi \left( x\right) \leq 1$ and $\varphi \left(
x\right) =0$\ for $\left| x\right| >2.$ Define $\varphi
_{\varepsilon }\left( x\right) =\varphi (\frac{x}{\varepsilon }),$
then $\varphi \left( x\right) =1$ for $\left| x\right| \leq
\varepsilon $\ and $\varphi _{\varepsilon }\left( x\right) =0$ for
$\left| x\right| >2\varepsilon $, so if $0<\rho <\varepsilon $,\
we will have
\begin{equation*}
u=\varphi _{\varepsilon }u\text{ \ , }\forall u\in C_{0}^{\infty }\left(
\Omega _{\rho }\right)
\end{equation*}
Let $\partial _{j}^{k}$\ denotes the derivation of order $k$\ with respect
to the variable $\xi _{j}$. It is well-known from the theory of
pseudodifferential operators that
\begin{equation}
P\left( ix_{j}u\right) =ix_{j}Pu+\left( \partial _{j}P\right) u\text{ ,}
\label{2}
\end{equation}
so, $\forall u\in C_{0}^{\infty }\left( \Omega _{\rho }\right) $, with $%
0<\rho <\varepsilon $, we have
\begin{equation}
P\left( ix_{j}u\right) =ix_{j}\varphi _{\varepsilon }\left( x\right)
Pu+\left( \partial _{j}P\right) u+T_{1}u\text{ , }  \label{3}
\end{equation}
where
\begin{equation*}
T_{1}=ix_{j}[P,\varphi _{\varepsilon }]
\end{equation*}
Then \ \
\begin{equation*}
\left\| \left( \partial _{j}P\right) u\right\| _{s}^{2}=\left( P\left(
ix_{j}u\right) ,\left( \partial _{j}P\right) u\right) _{s}-\left(
ix_{j}\varphi _{\varepsilon }\left( x\right) Pu,\left( \partial _{j}P\right)
u\right) _{s}-\left( T_{1}u,\left( \partial _{j}P\right) u\right) _{s}
\end{equation*}
It is easy to see that
\begin{equation*}
\left( P\left( ix_{j}u\right) ,\left( \partial _{j}P\right) u\right)
_{s}=\left( \overline{\left( \partial _{j}P\right) }\left( ix_{j}u\right) ,%
\overline{P}u\right) _{s}\text{ ,}
\end{equation*}
and consequently, we obtain
\begin{equation*}
\left\| \left( \partial _{j}P\right) u\right\| _{s}^{2}=\left( \overline{%
\left( \partial _{j}P\right) }\left( ix_{j}u\right) ,\overline{P}u\right)
_{s}-\left( ix_{j}\varphi _{\varepsilon }\left( x\right) Pu,\left( \partial
_{j}P\right) u\right) _{s}-\left( T_{1}u,\left( \partial _{j}P\right)
u\right) _{s}
\end{equation*}
From (\ref{3}), we have
\begin{equation}
\overline{\left( \partial _{j}P\right) }\left( ix_{j}u\right) =ix_{j}\varphi
_{\varepsilon }\left( x\right) \overline{\left( \partial _{j}P\right) }u+%
\overline{\left( \partial _{j}^{2}P\right) }u+T_{2}u\text{ \ ,}  \label{4}
\end{equation}
where
\begin{equation*}
T_{2}=ix_{j}[\overline{\left( \partial _{j}P\right) },\varphi _{\varepsilon
}]
\end{equation*}
Consequently, we have the following inequality
\begin{eqnarray}
\left\| \left( \partial _{j}P\right) u\right\| _{s}^{2} &\leq
&\left\| ix_{j}\varphi _{\varepsilon }\left( x\right)
\overline{\left( \partial _{j}P\right) }u\right\| _{s}\left\|
Pu\right\| _{s}+\left\| \left( \partial_{j}^{2}P\right) u\right\| _{s}\left\| Pu\right\| _{s}  \notag  \\
&&+\left\| ix_{j}\varphi _{\varepsilon }\left( x\right) Pu\right\|
_{s}\left\| \left( \partial _{j}P\right) u\right\| _{s}+\left| \left( T_{2}u,%
\overline{P}u\right) _{s}\right| +\left| \left( T_{1}u,\left(
\partial _{j}P\right) u\right) _{s}\right|  \label{5}
\end{eqnarray}
The lemma 2 gives
\begin{eqnarray*}
\left\| ix_{j}\varphi _{\varepsilon }\left( x\right) \overline{\left(
\partial _{j}P\right) }u\right\| _{s} &\leq &\max_{x}\left| ix_{j}\varphi
_{\varepsilon }\left( x\right) \right| \left\| \overline{\left( \partial
_{j}P\right) }u\right\| _{s}+c_{s,\sigma }\left( \varepsilon \right) \left\|
\overline{\left( \partial _{j}P\right) }u\right\| _{\sigma } \\
&\leq &2\varepsilon \left\| \left( \partial _{j}P\right) u\right\|
_{s}+c_{s,\sigma }\left( \varepsilon \right) \left\| u\right\| _{\sigma +m-1}%
\text{ , }\sigma <s-1\text{ \ ,}
\end{eqnarray*}
and
\begin{eqnarray*}
\left\| ix_{j}\varphi _{\varepsilon }\left( x\right) Pu\right\| _{s} &\leq
&\max_{x}\left| ix_{j}\varphi _{\varepsilon }\left( x\right) \right| \left\|
Pu\right\| _{s}+c_{s,\sigma }^{\prime }\left( \varepsilon \right) \left\|
Pu\right\| _{\sigma } \\
&\leq &2\varepsilon \left\| Pu\right\| _{s}+c_{s,\sigma }^{\prime }\left(
\varepsilon \right) \left\| u\right\| _{\sigma +m}\text{ , }\sigma <s-1
\end{eqnarray*}
For every real $t,$\ we have
\begin{equation*}
\left| \left( T_{2}u,\overline{P}u\right) _{s}\right| =\left| \left(
PT_{2}u,u\right) _{s}\right| \leq \left\| T_{2}u\right\| _{s-t+m}\left\|
u\right\| _{s+t}\text{ \ ,}
\end{equation*}
and
\begin{equation*}
\left| \left( T_{1}u,\left( \partial _{j}P\right) u\right) _{s}\right|
=\left| \left( \left( \overline{\partial _{j}P}\right) T_{1}u,u\right)
_{s}\right| \leq \left\| T_{1}u\right\| _{s-t+m-1}\left\| u\right\| _{s+t}
\end{equation*}
The above inequalities are resumed in the following one
\begin{eqnarray}
\left\| \left( \partial _{j}P\right) u\right\| _{s}^{2} &\leq
&4\varepsilon \left\| Pu\right\| _{s}\left\| \left( \partial
_{j}P\right) u\right\| _{s}+\left\| \left( \partial
_{j}^{2}P\right) u\right\| _{s}\left\| Pu\right\| _{s}+c_{s,\sigma
}^{\prime }\left( \varepsilon \right) \left\| \left( \partial
_{j}P\right) u\right\| _{s}\left\| u\right\| _{\sigma +m}+ \notag
\\
 &&+c_{s,\sigma }\left( \varepsilon \right) \left\|
Pu\right\| _{s}\left\| u\right\| _{\sigma +m-1}+\left\|
T_{1}u\right\| _{s-\tau +m-1}\left\| u\right\| _{s+\tau }+\left\|
T_{2}u\right\| _{s-t+m}\left\| u\right\| _{s+t}
\notag \\
&\leq &6\varepsilon \left\| Pu\right\| _{s}^{2}+4\varepsilon
\left\| \left(
\partial _{j}P\right) u\right\| _{s}^{2}+\frac{1}{8\varepsilon }\left\|
\left( \partial _{j}^{2}P\right) u\right\| _{s}^{2}+\frac{\left[
c_{s,\sigma }\left( \varepsilon \right) \right] ^{2}}{8\varepsilon
}\left\| u\right\|_{\sigma +m-1}^{2}+ \notag \\
&&+\frac{\left[ c_{s,\sigma }^{\prime }\left( \varepsilon \right)
\right] ^{2}}{8\varepsilon }\left\| u\right\| _{\sigma +m}^{2}+
\left\| u\right\|_{s+t}^{2}+\frac{1}{2}\left\|
 T_{1}u\right\| _{s-t+m-1}^{2}+\frac{1}{2}\left\| T_{2}u\right\| _{s-t+m}^{2} \notag
\\
\label{6}
\end{eqnarray}
Due to the lemma 5, the operators $T_{1}$\ et $T_{2}$\ are of
infinite orders. Let $\varepsilon >0$\ with $1-4\varepsilon >0$
and let $\sigma =s-\theta ,\ \theta >1,$ then there exists a
constant $c_{s,\theta }\left( \varepsilon \right) >0$\ such that
\begin{equation*}
\left\| \left( \partial _{j}P\right) u\right\| _{s}^{2}\leq \frac{%
6\varepsilon }{1-4\varepsilon }\left\| Pu\right\| _{s}^{2}+\frac{1}{%
(1-4\varepsilon )8\varepsilon }\left\| \left( \partial
_{j}^{2}P\right) u\right\| _{s}^{2}+c_{s,\theta }\left(
\varepsilon \right) \left\| u\right\| _{s+m-\theta }^{2}
\end{equation*}
Let $\delta >0$ and take $\varepsilon =\dfrac{\delta }{2\left( 2\delta
+3\right) }$ , then $\forall s\in \mathbb{R},\forall \theta >1,\forall
\delta >0,$ there exist $c_{1}\left( \delta \right) >0$ and $c_{1,s,\theta
}\left( \delta \right) >0,$
\begin{equation}
\left\| \left( \partial _{j}P\right) u\right\| _{s}^{2}\leq \delta \left\|
Pu\right\| _{s}^{2}+c_{1}\left( \delta \right) \left\| \left( \partial
_{j}^{2}P\right) u\right\| _{s}^{2}+c_{1,s,\theta }\left( \delta \right)
\left\| u\right\| _{s+m-\theta }^{2}  \label{7}
\end{equation}
$\forall u\in C_{0}^{\infty }\left( \Omega _{\rho }\right) ,\rho
<\varepsilon \leq \varepsilon _{1}\left( \delta \right) =\dfrac{\delta }{%
2\left( 2\delta +3\right) }$ .\smallskip \smallskip \smallskip

Let us show by induction that $\forall k\geq 1,\forall s\in \mathbb{R}%
,\forall \theta >1,\forall \delta >0,$ there exist $c_{k}\left( \delta
\right) >0,c_{k,s,\theta }\left( \delta \right) >0$ and $\varepsilon
_{k}\left( \delta \right) >0,\forall u\in C_{0}^{\infty }\left( \Omega
_{\rho }\right) ,\rho <\varepsilon \leq \varepsilon _{k}\left( \delta
\right) ,$ we have
\begin{equation}
\left\| \left( \partial _{j}^{k}P\right) u\right\| _{s}^{2}\leq \delta
\left\| Pu\right\| _{s}^{2}+c_{k}\left( \delta \right) \left\| \left(
\partial _{j}^{k+1}P\right) u\right\| _{s}^{2}+c_{k,s,\theta }\left( \delta
\right) \left\| u\right\| _{s+m-\theta }^{2}  \label{8}
\end{equation}
The case $k=1$ is true by (\ref{7}). Assume $\forall l\leq k-1,\forall s\in
\mathbb{R},\forall \theta >1,\forall \delta _{l}>0,$ there exist $%
c_{l}\left( \delta _{l}\right) >0,c_{l,s,\theta }\left( \delta _{l}\right)
>0,\forall u\in C_{0}^{\infty }\left( \Omega _{\rho }\right) ,\rho
<\varepsilon \leq \varepsilon _{l}\left( \delta _{l}\right) ,$ we have
\begin{equation}
\left\| \left( \partial _{j}^{l}P\right) u\right\| _{s}^{2}\leq \delta
_{l}\left\| Pu\right\| _{s}^{2}+c_{l}\left( \delta _{l}\right) \left\|
\left( \partial _{j}^{l+1}P\right) u\right\| _{s}^{2}+c_{l,s,\theta }\left(
\delta _{l}\right) \left\| u\right\| _{s+m-\theta }^{2}  \label{9}
\end{equation}
Apply the inequality (\ref{7}) to the operator $\left( \partial
_{j}^{k-1}P\right) $, i.e. $\forall $ $\delta ^{\prime }>0,$ there exist
positive constants $c_{1}\left( \delta ^{\prime }\right) ,c_{1,s,\theta
}\left( \delta ^{\prime }\right) $ and $\varepsilon _{1}(\delta ^{\prime
}),\forall u\in C_{0}^{\infty }\left( \Omega _{\rho }\right) ,\rho
<\varepsilon \leq \varepsilon _{1}\left( \delta ^{\prime }\right) ,$ we have
\
\begin{equation}
\left\| \left( \partial _{j}^{k}P\right) u\right\| _{s}^{2}\leq \delta
^{\prime }\left\| \left( \partial _{j}^{k-1}P\right) u\right\|
_{s}^{2}+c_{1}\left( \delta ^{\prime }\right) \left\| \left( \partial
_{j}^{k+1}P\right) u\right\| _{s}^{2}+c_{1,s,\theta }\left( \delta ^{\prime
}\right) \left\| u\right\| _{s+m-\theta }^{2}  \label{10}
\end{equation}
In (\ref{10}), we estimate $\left\| \left( \partial _{j}^{k-1}P\right)
u\right\| _{s}^{2}$ by the inequality (\ref{9}) with $l=k-1,$ then $\forall
u\in C_{0}^{\infty }\left( \Omega _{\rho }\right) ,\rho <\varepsilon \leq
\min \left\{ \varepsilon _{1}\left( \delta ^{\prime }\right) ,\varepsilon
_{k-1}\left( \delta _{k-1}\right) \right\} ,$ we obtain
\begin{eqnarray*}
\left\| \left( \partial _{j}^{k}P\right) u\right\| _{s}^{2} &\leq &\delta
^{\prime }\delta _{k-1}\left\| Pu\right\| _{s}^{2}+\delta ^{\prime
}c_{k-1}\left( \delta _{k-1}\right) \left\| \left( \partial _{j}^{k}P\right)
u\right\| _{s}^{2}+ \\
&&+\delta ^{\prime }c_{k-1,s,\theta }\left( \delta _{k-1}\right) \left\|
u\right\| _{s+m-\theta }^{2}+c_{1}\left( \delta ^{\prime }\right) \left\|
\left( \partial _{j}^{k+1}P\right) u\right\| _{s}^{2}+c_{1,s}\left( \delta
^{\prime }\right) \left\| u\right\| _{s+m-\theta }^{2}
\end{eqnarray*}
Choose $\delta ^{\prime }<\dfrac{1}{c_{k-1}\left( \delta _{k-1}\right) }$ ,
so $\forall u\in C_{0}^{\infty }\left( \Omega _{\rho }\right) ,\rho
<\varepsilon \leq \min \left\{ \varepsilon _{1}\left( \delta ^{\prime
}\right) ,\varepsilon _{k-1}\left( \delta _{k-1}\right) \right\} ,$ we have
\begin{eqnarray*}
\left\| \left( \partial _{j}^{k}P\right) u\right\| _{s}^{2} &\leq &\frac{%
\delta ^{\prime }\delta _{k-1}}{1-\delta ^{\prime }c_{k-1}\left( \delta
_{k-1}\right) }\left\| Pu\right\| _{s}^{2}+\frac{c_{1}\left( \delta ^{\prime
}\right) }{1-\delta ^{\prime }c_{k-1}\left( \delta _{k-1}\right) }\left\|
\left( \partial _{j}^{k+1}P\right) u\right\| _{s}^{2}+ \\
&&+\frac{c_{1,s}\left( \delta ^{\prime }\right) +\delta ^{\prime
}c_{k-1,s,\theta }\left( \delta _{k-1}\right) }{1-\delta ^{\prime
}c_{k-1}\left( \delta _{k-1}\right) }\left\| u\right\| _{s+m-\theta }^{2}
\end{eqnarray*}
Let $\delta >0$, take $\delta ^{\prime }=\dfrac{\delta }{\delta
_{k-1}+\delta c_{k-1}\left( \delta _{k-1}\right) }$ \ and
\begin{equation*}
c_{k}\left( \delta _{k}\right) =\dfrac{c_{1}\left( \delta ^{\prime }\right)
}{1-\delta ^{\prime }c_{k-1}\left( \delta _{k-1}\right) }\text{ \ ,}
\end{equation*}
and
\begin{equation*}
c_{k,s}\left( \delta _{k}\right) =\frac{c_{1,s}\left( \delta ^{\prime
}\right) +\delta ^{\prime }c_{k-1,s}\left( \delta _{k-1}\right) }{1-\delta
^{\prime }c_{k-1}\left( \delta _{k-1}\right) }
\end{equation*}
Then, we obtain
\begin{equation*}
\left\| \left( \partial _{j}^{k}P\right) u\right\| _{s}^{2}\leq \delta
_{k}\left\| Pu\right\| _{s}^{2}+c_{k}\left( \delta _{k}\right) \left\|
\left( \partial _{j}^{k+1}P\right) u\right\| _{s}^{2}+c_{k,s,\theta }\left(
\delta _{k}\right) \left\| u\right\| _{s+m-\theta }^{2}
\end{equation*}
$\forall u\in C_{0}^{\infty }\left( \Omega _{\rho }\right) ,\rho
<\varepsilon \leq \varepsilon _{k}\left( \delta _{k}\right) =\min \left\{
\varepsilon _{1}\left( \dfrac{\delta }{\delta _{k-1}+\delta c_{k-1}\left(
\delta _{k-1}\right) }\right) ,\varepsilon _{k-1}\left( \delta _{k}\right)
\right\} $. We have proved the inequality (\ref{9}) for $l=k.$ So the
estimate (\ref{8}) is true.

Let $\delta _{k},c_{k}\left( \delta _{k}\right) ,c_{k,s,\theta }\left(
\delta _{k}\right) $ et$\ \varepsilon _{k}\left( \delta _{k}\right)
,k=1,2,...,l$ , be the respective constants of the right member of the
estimates (\ref{8}).$\ $Iterating these inequalities, then $\forall u\in
C_{0}^{\infty }\left( \Omega _{\rho }\right) ,\rho <\varepsilon \leq \min
\left\{ \varepsilon _{1}\left( \delta _{1}\right) ,...,\varepsilon
_{l}\left( \delta _{l}\right) \right\} ,$ $\forall \ l\geq 1$,\ we obtain
the following one,
\begin{eqnarray}
\left\| \left( \partial _{j}P\right) u\right\| _{s}^{2} &\leq &\left( \delta
_{1}+c_{1}\left( \delta _{1}\right) \delta _{2}+....+c_{1}\left( \delta
_{1}\right) c_{2}\left( \delta _{2}\right) ...c_{l-1}\left( \delta
_{l-1}\right) \delta _{l}\right) \left\| Pu\right\| _{s}^{2}+  \notag \\
&&+c_{1}\left( \delta _{1}\right) c_{2}\left( \delta _{2}\right)
....c_{l}\left( \delta _{l}\right) \left\| \left( \partial
_{j}^{l+1}P\right) u\right\| _{s}^{2}+  \notag \\
&&\left( c_{1,s,\theta }\left( \delta _{1}\right) +c_{2,s,\theta
}\left( \delta _{2}\right) c_{1}\left( \delta _{1}\right)
+c_{l,s,\theta }\left( \delta _{l}\right) c_{1}\left( \delta
_{1}\right) ...c_{l-1}\left( \delta _{l-1}\right) \right) \left\|
u\right\| _{s+m-\theta }^{2}  \notag \\
\label{11}
\end{eqnarray}
Let $\delta >0$, choose $\delta _{1},...,\delta _{l}$ from the following
equations \ \ \
\begin{equation*}
\delta _{1}=\frac{\delta }{l},c_{1}\left( \delta _{1}\right) \delta _{2}=%
\frac{\delta }{l},...,c_{1}\left( \delta _{1}\right) c_{2}\left( \delta
_{2}\right) ...c_{l-1}\left( \delta _{l-1}\right) \delta _{l}=\frac{\delta }{%
l}\text{ },
\end{equation*}
and define the constants $c_{l}\left( \delta \right) $ and $c_{l,s,\theta
}\left( \delta \right) $ respectively as the coefficients of the terms $%
\left\| \left( \partial _{j}^{t+1}P\right) u\right\| _{s}^{2}$ and $\left\|
u\right\| _{s+m-\theta }^{2}$ in the inequality (\ref{11}), then $\forall
l\geq 1,\forall s\in \mathbb{R},\forall \theta >1,\forall \delta >0,$ there
exist\ $c_{l}\left( \delta \right) >0,c_{l,s,\theta }\left( \delta \right)
>0 $ and $\varepsilon _{l}\left( \delta \right) >0,$ $\forall u\in
C_{0}^{\infty }\left( \Omega _{\rho }\right) ,\rho <\varepsilon \leq
\varepsilon _{l}\left( \delta \right) ,$ we have
\begin{equation}
\left\| \left( \partial _{j}P\right) u\right\| _{s}^{2}\leq \delta \left\|
Pu\right\| _{s}^{2}+c_{l}\left( \delta \right) \left\| \left( \partial
_{j}^{l+1}P\right) u\right\| _{s}^{2}+c_{l,s,\theta }\left( \delta \right)
\left\| u\right\| _{s+m-\theta }^{2}  \label{12}
\end{equation}
Choose $l\in \mathbb{N}$ with $l\geq \theta -1>0,$ then there is $%
c_{s,\theta }\left( \delta \right) >0,\forall u\in C_{0}^{\infty }\left(
\Omega _{\rho }\right) ,\rho <\varepsilon \leq \min \left\{ \varepsilon
_{1}\left( \delta _{1}\right) ,...,\varepsilon _{l}\left( \delta _{l}\right)
\right\} ,$
\begin{equation}
\left\| \left( \partial _{j}P\right) u\right\| _{s}^{2}\leq \delta \left\|
Pu\right\| _{s}^{2}+c_{s,\theta }\left( \delta \right) \left\| u\right\|
_{s+m-\theta }^{2}  \label{13}
\end{equation}
The inequality (\ref{13}) is true for $\theta =1,$\ because the operator $%
\left( \partial _{j}P\right) $\ is of order $s+m-1.$ Finally, we have proved
that $\forall s\in \mathbb{R},\forall \theta \geq 1,\forall \delta
>0,\exists \rho >0,\exists c>0,\forall u\in C_{0}^{\infty }\left( \Omega
\right) ,diam(\Omega )<\rho ,$\
\begin{equation}
\left\| \left( \partial _{j}P\right) u\right\| _{s}^{2}\leq \delta \left\|
Pu\right\| _{s}^{2}+c\left\| u\right\| _{s+m-\theta }^{2}  \label{14}
\end{equation}

Let $\alpha =\left( \alpha _{1},...,\alpha _{n}\right) $ and
$\alpha ^{\prime }=\left( \alpha _{1},...,\alpha _{j-1},\alpha
_{j}+1,\alpha _{j+1}...,\alpha _{n}\right) $ be given
multi-indices. Assume as an hypothesis of induction : $\forall
s\in \mathbb{R},\forall \theta \geq 1,\forall \delta >0,\forall
\alpha \in \mathbb{Z}_{+}^{n},\exists \rho
>0,\exists c>0,\forall u\in C_{0}^{\infty }\left( \Omega \right)
,diam(\Omega )<\rho ,$
\begin{equation}
\left\| \left( \partial _{j}^{\alpha }P\right) u\right\| _{s}^{2}\leq \delta
\left\| Pu\right\| _{s}^{2}+c\left\| u\right\| _{s+m-\theta }^{2}\text{ \ ,}
\label{15}
\end{equation}
is true. Apply the inequality (\ref{14}) to the operator $\left( \partial
^{\alpha }P\right) ,$ then we have
\begin{equation*}
\left\| \left( \partial _{j}^{\alpha ^{\prime }}P\right) u\right\|
_{s}^{2}\leq \delta ^{\prime }\left\| \left( \partial ^{\alpha }P\right)
u\right\| _{s}^{2}+c^{\prime }\left\| u\right\| _{s+m-\theta }^{2}\text{ , }%
u\in C_{0}^{\infty }\left( \Omega ^{\prime }\right) \text{ \ ,}
\end{equation*}
where $\Omega ^{\prime }$ depends on $\delta ^{\prime }.$ From the
hypothesis of induction for the operator$\left( \partial ^{\alpha }P\right) $%
, we obtain that for every $\delta >0,$ there is $\rho >0$ such that
\begin{equation*}
\left\| \left( \partial _{j}^{\alpha ^{\prime }}P\right) u\right\|
_{s}^{2}\leq \delta ^{\prime }\delta \left\| Pu\right\| _{s}^{2}+\delta
^{\prime }c\left\| u\right\| _{s+m-\theta }^{2}+c^{\prime }\left\| u\right\|
_{s+m-\theta }^{2}\text{ \ ,}
\end{equation*}
$u\in C_{0}^{\infty }\left( \Omega \cap \Omega ^{\prime }\right) ,diam\left(
\Omega \right) <\rho .$ Let $\gamma >0,$ choose $\delta ^{\prime }=\dfrac{%
\gamma }{\delta }$ , we obtain then the inequality (\ref{15})\ for
$\alpha ^{\prime }$. This ends the proof of the theorem
\end{proof}

\begin{remark}
As the operator $\left( \partial ^{\alpha }P\right) $\ is of order $m-\left|
\alpha \right| ,$ the estimate (\ref{1}) is trivial for$\ \theta \leq \left|
\alpha \right| ,$ the theorem is then restated in the following way :\newline
Let $P\left( D\right) $ be a pseudodifferential operator of the class $%
S^{m},s^{\prime },s\in \mathbb{R},$ and$\ \alpha \in \mathbb{N}^{n}$, $%
s^{\prime }<s+m-\left| \alpha \right|,$ then $\forall \delta >0,\exists \rho
>0,\exists c>0,\forall u\in C_{0}^{\infty }\left( \Omega \right)
,diam(\Omega )<\rho ,$
\begin{equation*}
\left\| \left( \partial ^{\alpha }P\right) \left( D\right) u\right\|
_{s}\leq \delta \left\| P\left( D\right) u\right\| _{s}+c\left\| u\right\|
_{s^{\prime }}\text{ \ ,}
\end{equation*}
In this form, the theorem reminds the so called Erhling's inequality.
\end{remark}


\begin{thebibliography}{9}
\bibitem{B-K}  Bouzar C., Kuleshov A. A., Local solvability of
pseudodifferential equations of constant strength I. Differential Equations
24, no. 5, 548--553, (1988).

\bibitem{B}  Bouzar C., Local estimates for pseudodifferential operators.
Dokl. Nats. Akad. Nauk Belarusi 44, no. 4, 18--20, (2000)

\bibitem{H1}  H\"{o}rmander L., On the theory of general partial
differential operators. Acta Math., 94, 161-248 (1955).

\bibitem{H2}  H\"{o}rmander L., The analysis of linear partial differential
operators I. \ Springer, Second Edition 1990.

\bibitem{P}  Paneah B., Pseudodifferential operators of constant strength in
the main. Math. Sbornik 2:2, 179-198, (1967).

\bibitem{Pe}  Peetre J., A proof of the hypoellipticity of formally
hypoelliptic differential operators. Comm. Pure Appl. Math., 14, 737-744,
(1961).

\bibitem{T}  Treves F., Introduction to pseudodifferential and Fourier
integral operatoprs, Volume I, Plenum Press, 1982.
\end{thebibliography}
\end{document}